\documentclass[preprint,12pt]{elsarticle}
\usepackage{amsmath}
\usepackage{amsmath,amssymb,amsthm}

\newcommand{\re}{\mathbb{R}}
\newcommand{\co}{\mathbb{C}}

\newcommand{\rp}{\mbox{Re}}
\newcommand{\ip}{\mbox{Im}}

\newtheorem{thm}{Theorem}[section]

\newtheorem{lem}[thm]{Lemma}
\newtheorem{cor}[thm]{Corollary}
\newtheorem{rem}[thm]{Remark}

\theoremstyle{definition}

\usepackage{amssymb}

\journal{}

\begin{document}

\begin{frontmatter}



\title{A boundary Schwarz lemma for holomorphic mappings from the polydisc to the unit ball}

\author[1,3]{Yang Liu\fnref{label2}}
\address[1]{Department of Mathematics, Zhejiang Normal
University, Jinhua 321004, China} \fntext[label2]{Email: liuyang@zjnu.edu.cn}

\author[2]{Zhihua Chen\fnref{label3}}
\address[2]{Department of Mathematics, Tongji University, Shanghai
200092, China}\fntext[label3]{Email: zzzhhc@tongji.edu.cn}

\author[3,4]{Yifei Pan\fnref{label3}}
\address[3]{Department of Mathematical Sciences,
Indiana University-Purdue University Fort Wayne,
Fort Wayne, IN 46805-1499, USA}
\address[4]{School of Mathematics and Informatics, Jiangxi Normal University, Nanchang 330022, China}
\fntext[label4]{Email: pan@ipfw.edu}

\address{}

\begin{abstract}
In this paper, we prove a general Schwarz lemma at the boundary for holomorphic mappings from the polydisc to the unit ball in any dimensions. For the special case of one complex variable, the obtained results give the classic boundary Schwarz lemma.
\end{abstract}


\end{frontmatter}


\section{ Introduction}\label{section1}

The Schwarz lemma is known as one of the most important results in complex analysis. The Schwarz-Pick lemma is a variant of the Schwarz lemma, which essentially states that a holomorphic self-map of the unit disk decreases the distance of points in the Poincar\'e metric. It has been generalized to the derivatives of arbitrary order in one complex variable \cite{ruscheweyh1985two,dai2008note}. Concerning several complex variables, Rudin \cite{rudin1969function} gave a first derivative estimate for the bounded holomorphic functions on the polydisc, which is really a precursor to Schwarz-Pick estimate in high dimensions. \cite{knese2007schwarz} studied the extremal problem for
the holomorphic functions from the polydisk to the disk. later, \cite{chen2010schwarz,DCP1,sai2010,liu2011schwarz} generalized the result of \cite{ruscheweyh1985two,dai2008note} to the holomorphic mappings on the unit ball and polydisc in $\co^n$.

Schwarz lemma at the boundary is another active topic in complex analysis \cite{garnett1981bounded}, which was applied to geometric function theory of one complex variable and several complex
variables \cite{krantz2011schwarz,lwt}. The following result is the classic boundary version of Schwarz lemma in one complex variable.
\begin{thm}[\cite{garnett1981bounded}]\label{thmg}
Let $D$ be the unit disk in $\co$, and let $f$ be the self-holomorphic mapping of $D$. If $f$ is holomorphic at $z=1$ with $f(0)=0$ and $f(1)=1$, then $f'(1)\geq 1$.
\end{thm}
In \cite{burns1994rigidity}, Burns and Krantz gave a new Schwarz lemma at the
boundary and obtained rigidity result for holomorphic mappings. In \cite{baracco2006burns}, the Burns-Krantz type theorem was extended to mappings defined only on one side as germs at a boundary point and not necessarily sending any domain to itself. There are also some other multidimensional generalization of the Schwarz lemma at the boundary, see \cite{huang1995boundary,lwt}.

It is a well-known result that there are no biholomorphic mappings between $D^n$ and $B^n$. Therefore, holomorphic mappings
between $D^n$ and $B^n$ are interesting in several complex variables. There have been many results proved about such mappings
since the 1970s, see \cite{db2,db1,stessin2006composition,liu2011schwarz}. Before giving the proof of the boundary version of Schwarz lemma for holomorphic mappings
between $D^n$ and $B^N$, we give some notations and definitions first.

For any $z=(z_1,...,z_n)^T,~w=(w_1,...,w_n)^T\in \co^n$, the inner product and the corresponding norm are given by $\langle z,~w\rangle=\sum_{j=1}^{n}z_j\overline{w_j},~||z||=\langle z,z\rangle^{\frac{1}{2}}$. In addition, denote $||z||_\infty=\max_{1\leq j\leq n}|z_j|$.  Let $D^n\in \co^n$ be the polydisc in the $n$-dimensional complex space, which is described by $D^n=\left\{z\in\co^n|||z||_\infty<1\right\}$. $\partial D^n=\left\{z\in\co^n|||z||_\infty=1\right\}$ denotes the topological boundary of $D^n$, and the distinguished boundary of $D^n$ is given by $T^n=\left\{z\in\co^n||z_j|=1,~1\leq j\leq n\right\}$, see \cite{rudin1969function}. Denote $H(D^n,B^N)$ by the set of all holomorphic mappings from $D^n$ to $B^N$. For any $f\in H(D^n,B^N)$, we denote it by $f=(f_1,...,f_N)^T$ and the Jacobian matrix of $f$ at $z\in D^n$ is given by $$J_f(z)=\left(\frac{\partial f_i}{\partial z_j}\right)_{N\times n}.$$
For a bounded domain $V\in \co^n$, $C^{\alpha}(V)$ for $0<\alpha<1$ is the set of all functions $f$ on $V$ for which
$$\sup\left\{{\frac{|f(z)-f(z')|}{|z-z'|^{\alpha}} \bigg| z,z' \in \bar V} \right\}$$
is finite.
$C^{k+\alpha}(V)$ is the set of all functions $f$ on $V$ whose $k^{\textup{th}}$ order partial derivatives exist and belong to $C^{\alpha}(V)$ for an integer $k \geq 0$.

\begin{rem}
From the definition of the norm on $D^n$, the boundary points  could be classified into several kinds. Let $z_0=(z_1,...,z_n)^T\in \partial D^n$. If there are only $r$ components of $z_0$ whose norm equals to $1$, then we denote the set of all this kind of boundary points by $E_r$, $1\leq r\leq n$. It is easy to see that there are $n$ different $E_r$ for $\partial D^n$, i.e., $$\bigcup_{1\leq r\leq n}E_r=\partial D^n.$$ Especially, $E_n=T^n$ which is the distinguished boundary of $D^n$.

\end{rem}

A natural question for
Schwarz lemma problem and rigidity problem at the boundary was asked by Krantz in \cite{krantz2011schwarz} for non-equidimensional mappings. In this paper, we study the mapping $f\in H(D^n,B^N)$ for any $n,N\geq 1$. Our main results are listed as follows.
\begin{thm}\label{thm1}
Let $f\in H(D^n,B^N)$ for any $n,N\geq 1$. Given $z_0\in \partial D^n$. Assume $z_0\in E_r$ with the first $r$ components at the boundary of $D$ for some $1\leq r\leq n$. If $f$ is $C^{1+\alpha}$ at $z_0$ and $f(z_0)=w_0\in \partial B^N$, then there exist a sequence of nonnegative real numbers $\gamma_1,...,\gamma_r$ satisfying $\sum_{j=1}^{r}\gamma_j\geq 1$ and $\lambda\in \re$ such that
\begin{equation}\label{cc}
\begin{split}
\overline{J_f(z_0)}^Tw_0=\lambda \mathrm {diag}(\gamma_1,...,\gamma_r,0,...,0) z_0,
\end{split}
\end{equation}
where $\lambda=\frac{|1-\bar{a}^Tw_0|^2}{1-||a||^2}>0$, $ a=f(0)$ and $\mathrm {diag}$ represents the diagonal matrix.
\end{thm}

\begin{rem}\label{}
Especially, when $r=n$, $z_0\in T^n$. Denote $e^n_i$ by the $i$-th column of identity matrix $I_n$ with degree $n$. If $z_0\in E_r$ is given by other expression such as $z_0=\sum_{j=1}^{r}e^n_{l_j}$ where $l_j$ are different to each other and $i_j\in \{1,...,n\}$ for $1\leq j\leq r$. Then from the proof of the theorem, $\mathrm {diag}(\gamma_1,...,\gamma_r,0,...,0)$ should be replaced by the diagonal matrix $M$ with the $i_j$-th low and $i_j$-th column element $M(i_j,i_j)=\gamma_{i_j}$ ($j=1,...,r$), otherwise being $0$.

\end{rem}
For $r=1$, we have the following corollary.
\begin{cor}\label{cor1}
Let $f\in H(D^n,B^N)$ for any $n\geq 1$. Given $z_0\in E_1\subset \partial D^n$. If $f$ is $C^{1+\alpha}$ at $z_0$ and $f(z_0)=w_0\in \partial B^N$, then there exists a real number $\lambda\in \re$ such that
\begin{equation}\label{cc}
\begin{split}
\overline{J_f(z_0)}^Tw_0=\lambda z_0,
\end{split}
\end{equation}
where $\lambda\geq\frac{|1-\bar{a}^Tw_0|^2}{1-||a||^2}>0$ and $ a=f(0)$.
\end{cor}

For $N=1$, the following result is obtained.
\begin{cor}\label{cor}
Let $f\in H(D^n,D)$ for any $n\geq 1$. Given $z_0\in \partial D^n$. Assume $z_0\in E_r$ with the first $r$ components in the boundary of $D$ for some $1\leq r\leq n$. If $f$ is $C^{1+\alpha}$ at $z_0$ and $f(z_0)=e^{i\theta}\in \partial D$ for $0\leq \theta\leq 2\pi$. Then there exist a sequence of nonnegative real numbers $\gamma_1,...,\gamma_r$ satisfying $\sum_{j=1}^{r}\gamma_j\geq 1$ and $\lambda\in \re$ such that
\begin{equation}\label{cc}
\begin{split}
\overline{J_f(z_0)}^Te^{i\theta}=\lambda \mathrm {diag}(\gamma_1,...,\gamma_r,0,...,0) z_0,
\end{split}
\end{equation}
where $\lambda=\frac{|1-\bar{a}e^{i\theta}|^2}{1-|a|^2}>0$ and $ a=f(0)$.
\end{cor}
\begin{rem}\label{}
This theorem is a general Schwarz lemma at the boundary for holomorphic mappings from the polydisc to the unit ball in any dimensions. For $n=N=1$, Corollary \ref{cor1} gives Theorem \ref{thmg} in \cite{garnett1981bounded}. The smooth condition of $f$ is $C^{1+\alpha}$ at $z_0$ here.
\end{rem}



The Kobayashi distances for the polydisc in $\co^n$ is given as follows \cite{jarnicki1993invariant}.
\begin{equation}\label{01}
\begin{split}
K_{D^n}(z,w)=\frac{1}{2}\log \frac{1+||\varphi_z(w)||_\infty}{1-||\varphi_z(w)||_\infty},~z,w\in D^n
\end{split}
\end{equation}
where $\varphi_z(w)$ is the automorphism of $D^n$ given by
$$\varphi_z(w)=\left(\frac{w_1-z_1}{1-\bar z_1w_1},...,\frac{w_n-z_n}{1-\bar z_nw_n}\right)^T.$$
Meanwhile, the Kobayashi distances for the unit ball in $\co^N$ could be expressed by
\begin{equation}\label{02}
\begin{split}
K_{B^N}(z,w)=\frac{1}{2}\log \frac{1+||\phi_z(w)||}{1-||\phi_z(w)||}, ~z,w\in B^N
\end{split}
\end{equation}
where $\phi_z(w)$ is the automorphism of $B^N$, and
$$
\phi_z(w)=\frac{z-P_z(w)-s_zQ_z(w)}{1-\langle w,z\rangle},~z,w\in B^N
$$
with $P_z$ being the orthogonal projection of $\co^n$ by $$P_z(w)=\frac{\langle w,z\rangle}{||z||^2}z,~\mbox{if}~z\neq 0,$$
and $Q_z=I-P_z$, as well as $s_z=\sqrt{1-||z||^2}$, see \cite{rudin}.
It is found that $\phi_z(0)=z,~\phi_z(z)=0$ and $\phi_z=\phi_z^{-1}$.

It is the fact that the Kobayashi distance non-increases under holomorphic mappings \cite{kobayashi1967intrinsic}.
Consider the mapping $f\in H(D^n,B^N)$ and $f(0)=0$, then from (\ref{01}) and (\ref{02})
$K_{B^N}(f(0),f(w))\leq K_{D^n}(0,w),$ i.e.,
$$\frac{1}{2}\log \frac{1+||f(w)||}{1-||f(w)||}\leq \frac{1}{2}\log \frac{1+||w||_\infty}{1-||w||_\infty}.$$
Since $t\rightarrow \frac{1}{2}\log \frac{1+t}{1-t}$ is a increasing function for $t\in [0,1)$, we have $||f(w)||\leq ||w||_\infty$. Therefore,
 the following lemma is obtained which would play an important role in the proof of main results.

\begin{lem}\label{lem2}
If $f\in H(D^n,B^N)$ for any $n,N\geq 1$. If $f(0)=0$, then $||f(w)||\leq ||w||_\infty,~w\in D^n.$
\end{lem}
\section{Proof of Theorem \ref{thm1}}
\begin{proof} We will prove the theorem in the following five steps.

{\bf Step 1.}  Let $z_0\in E_r\subset \partial D^n$. Without loss of generality, we assume $z_0=\sum_{i=1}^{r}e_i^n$. $f$ is $C^{1+\alpha}$ in a neighborhood $V$ of $z_0$. Moreover, we assume $f(0)=0$ and $f(z_0)=w_0=e^N_1\in\partial B^N$.

Let $p=z_0$, $q_j=-\sum_{i=1}^{r}e_i^n+ike_j^n$ for $1\leq j\leq r$. Then
$p+tq_j=(1-t)z_0+ikte_j^n$ for $t\in \re$. $||p+tq_j||_\infty<1\Leftrightarrow |1-t+ikt|<1$ and $|1-t|<1\Leftrightarrow 0<t<\frac{2}{1+k^2}$, which means that for a given $k\in \re$ when $t \rightarrow 0^+$, $p+tq_j\in  D^n\cap V$. For such $t$, taking the Taylor expansion of $f\left(p+tq_j\right)$ at $t=0$, we have
$$f\left(p+tq_j\right)=w_0+J_f(z_0)q_jt+O(t^{1+\alpha}).$$
By Lemma \ref{lem2},
$$||f\left(p+tq_j\right)||^2=||w_0+J_f(z_0)q_jt+O(t^{1+\alpha})||^2\leq ||p+tq_j||_\infty^2,$$
i.e.,
$$1+2\rp\left(\overline{w_0}^TJ_f(z_0)q_jt\right)+O(t^{1+\alpha})\leq 1-2t+O(t^2).$$
Substitute $w_0=e^N_1$ and let $t \rightarrow 0^+$, we have
$$\rp\left(\overline{e^N_1}^TJ_f(z_0)(-\sum_{i=1}^{r}e_i^n+ike_j^n)\right)\leq -1,$$
i.e.,
$$\rp\left(-\sum_{i=1}^{r}\frac{\partial f_1(z_0)}{\partial z_i}+ik\frac{\partial f_1(z_0)}{\partial z_j}\right)\leq -1.$$
Let $\frac{\partial f_1(z_0)}{\partial z_j}=\rp \frac{\partial f_1(z_0)}{\partial z_j}+i\ip \frac{\partial f_1(z_0)}{\partial z_j}$, then from the above inequality, one gets
$$-\rp\sum_{i=1}^{r}\frac{\partial f_1(z_0)}{\partial z_i}-k\ip\frac{\partial f_1(z_0)}{\partial z_j}\leq -1,$$
i.e.,
\begin{equation}\label{1}
\begin{split}
-k\ip\frac{\partial f_1(z_0)}{\partial z_j}\leq \rp\sum_{i=1}^{r}\frac{\partial f_1(z_0)}{\partial z_i}-1.
\end{split}
\end{equation}
Since (\ref{1}) is valid for any $k\in \mathbb R$ so that
$$\ip\frac{\partial f_1(z_0)}{\partial z_j}=0,~1\leq j\leq r,$$
which gives
\begin{equation}\label{1.5}
\begin{split}\rp\sum_{i=1}^{r}\frac{\partial f_1(z_0)}{\partial z_i}\geq 1,\end{split}
\end{equation} and
\begin{equation}\label{2}
\begin{split}
\frac{\partial f_1(z_0)}{\partial z_j}=\rp\frac{\partial f_1(z_0)}{\partial z_j},~1\leq j\leq r.\end{split}
\end{equation}

{\bf Step 2.}
Let $p=z_0$, $q_j=-\sum_{i=1}^{r}e_i^n+ke_j^n$ for $1\leq j\leq r$ and $k\leq 0$. Then
$p+tq_j=(1-t)z_0+kte_j^n$ for $t\in \re$. $||p+tq_j||_\infty<1\Leftrightarrow |1-t+kt|<1$ and $|1-t|<1\Leftrightarrow 0<t<\frac{2}{1-k}$, which means that for any given $k\leq 0$, there are $t \rightarrow 0^+$ such that $p+tq_j\in  D^n\cap V$. Taking the Taylor expansion of $f\left(p+tq_j\right)$ at $t=0$, we have
$$||f\left(p+tq_j\right)||^2=||w_0+J_f(z_0)q_jt+O(t^{1+\alpha})||^2\leq ||p+tq_j||_\infty^2,$$
i.e.,
$$1+2\rp\left(\overline{w_0}^TJ_f(z_0)q_jt\right)+O(t^{1+\alpha})\leq 1-2t+t^2.$$
Substitute $w_0=e^N_1$ and let $t \rightarrow 0^+$, we have
$$\rp\left(\overline{e^N_1}^TJ_f(z_0)(-\sum_{i=1}^{r}e_i^n+ke_j^n)\right)\leq -1,$$
i.e.,
$$\rp\left(-\sum_{i=1}^{r}\frac{\partial f_1(z_0)}{\partial z_i}+k\frac{\partial f_1(z_0)}{\partial z_j}\right)\leq -1.$$
From (\ref{2}), one gets
$$-\sum_{i=1}^{r}\frac{\partial f_1(z_0)}{\partial z_i}+k\frac{\partial f_1(z_0)}{\partial z_j}\leq -1,$$
which equals to
$$k\frac{\partial f_1(z_0)}{\partial z_j}\leq \sum_{i=1}^{r}\frac{\partial f_1(z_0)}{\partial z_i}-1.$$
The right side of the above inequality is nonnegative from (\ref{1.5}) and (\ref{2}). Since $k\leq0$ is arbitrary, it is obtained that \begin{equation}\label{2.5}
\begin{split}
\frac{\partial f_1(z_0)}{\partial z_j}\geq 0,~1\leq j\leq r.\end{split}
\end{equation}

{\bf Step 3.}
Let $q_l=-\sum_{j=1}^{r}e^n_j+ike^n_l$ for $r+1\leq l\leq n$ and $k\neq 0\in \mathbb R$. Then
$p+tq_l=(1-t)\sum_{j=1}^{r}e^n_j+ikte^n_l$ for $t\in \re$. $||p+tq_l||_\infty<1\Leftrightarrow |1-t|<1$ and $|ikt|^2<1\Leftrightarrow 0<t<\min\{\frac{1}{k^2},2\}$. Therefore,  given a $k\neq 0\in \re$, there exist $t \rightarrow 0^+$ such that $p+tq\in D^n\cap V$.
Similarly, take the Taylor expansion of $f\left(p+tq_l\right)$ at $t=0$, we have
$$f\left(p+tq_l\right)=w_0+J_f(z_0)q_lt+O(t^{1+\alpha}).$$
By Lemma \ref{lem2},
$$||f\left(p+tq_l\right)||^2=||w_0+J_f(z_0)q_lt+O(t^{1+\alpha})||^2\leq ||p+tq_l||_\infty^2,$$
i.e.,
$$1+2\rp\left(\overline{w_0}^TJ_f(z_0)(-\sum_{j=1}^{r}e^n_j+ike^n_l)t\right)+O(t^{1+\alpha})\leq 1-2t+O(t^{2}).$$
Substitute $w_0=e^N_1$ and let $t \rightarrow 0^+$, we have
$$\rp\left(\overline{e^N_1}^TJ_f(z_0)(-\sum_{j=1}^{r}e^n_j+ike^n_l)\right)\leq -1,$$
i.e.,
$$\rp\left(-\sum_{j=1}^{r}\frac{\partial f_1(z_0)}{\partial z_j}+ik\frac{\partial f_1(z_0)}{\partial z_l}\right)\leq -1.$$
From the above inequality as well as inequality (\ref{2}), one has
$$k\ip\frac{\partial f_1(z_0)}{\partial z_l}\leq \sum_{j=1}^{r}\frac{\partial f_1(z_0)}{\partial z_j}-1.$$
Since the right side of above inequality is a nonnegative scalar, with the similar argument to Step 1, we have
$$\ip\frac{\partial f_1(z_0)}{\partial z_l}=0,~r+1\leq l\leq n.$$
Meanwhile, if we assume $p=z_0$, $q=-\sum_{j=1}^{r}e^n_j+ke^n_l$ for $r+1\leq l\leq n$ and any $k\neq 0\in \mathbb R$. It is easy to find
$$\rp\frac{\partial f_1(z_0)}{\partial z_l}=0,~r+1\leq l\leq n$$
as well.
Therefore,
\begin{equation}\label{3}
\begin{split}
\frac{\partial f_1(z_0)}{\partial z_l}=0,~r+1\leq l\leq n.
\end{split}
\end{equation}
As a result of (\ref{1.5}), (\ref{2.5}) and (\ref{3}), we have
\begin{equation}\label{4}
\begin{split}
\overline{J_f(z_0)}^Tw_0=\mathrm {diag}(\lambda_1,...,\lambda_r,0,...,0) z_0
\end{split}
\end{equation}
for $w_0=e^N_1,~z_0=\sum_{j=1}^{r}e^n_j$ and $\lambda_j=\frac{\partial f_1(z_0)}{\partial z_j}\geq 0$ with $\sum_{j=1}^{r}\lambda_j\geq 1$ and $1\leq j\leq r$.

{\bf Step 4.}
Now let $z_0\in E_r$ be any given point at $\partial D^n$ with the first $r$ components in the boundary of $D$, i.e., $z_0$ is not necessary $\sum_{j=1}^{r}e^n_j$. Then there exists a special kind of diagonal unitary matrix $U_{z_0}$ such that $U_{z_0}(z_0)=\sum_{j=1}^{r}e^n_j$. Assume $f(0)=0,~f(z_0)=w_0$ and $w_0$ is not necessary $e^N_1$ at $\partial B^N$. Similarly, there is a $U_{w_0}$ such that $U_{w_0}(w_0)=e^N_1$. Denote $$g(z)=U_{w_0}\circ f\circ \overline{U_{z_0}}^T,$$ then $g(0)=0,~g(\sum_{j=1}^{r}e^n_j)=e^N_1.$ Moreover,
\begin{equation}\label{5}
\begin{split}
J_g(z)=U_{w_0}J_f(\overline{U_{z_0}}^Tz)\overline{U_{z_0}}^T.
\end{split}
\end{equation}
 From Steps 2 and 3, we have
$$
\overline{J_g(\sum_{j=1}^{r}e^n_j)}^Te^N_1=\mathrm {diag}(\lambda_1,...,\lambda_r,0,...,0) \sum_{j=1}^{r}e^n_j
$$
where $\lambda_j=\frac{\partial g_1(z_0)}{\partial z_j}\geq 0$ with $\sum_{j=1}^{r}\lambda_j\geq 1$ and $1\leq j\leq r$, which equals to
$$
\overline{U_{w_0}J_f(\overline{U_{z_0}}^T\sum_{j=1}^{r}e^n_j)\overline{U_{z_0}}^T}^Te^N_1=\mathrm {diag}(\lambda_1,...,\lambda_r,0,...,0) \sum_{j=1}^{r}e^n_j,
$$
i.e.,
$$
U_{z_0}\overline{J_f(z_0)}^T\overline{U_{w_0}}^T e^N_1=\mathrm {diag}(\lambda_1,...,\lambda_r,0,...,0) \sum_{j=1}^{r}e^n_j.
$$
Multiplying $\overline{U_{z_0}}^T$ at both sides of the above equation gives
$$
\overline{U_{z_0}}^TU_{z_0}\overline{J_f(z_0)}^T\overline{U_{w_0}}^T e^N_1=\overline{U_{z_0}}^T\mathrm {diag}(\lambda_1,...,\lambda_r,0,...,0) \sum_{j=1}^{r}e^n_j.
$$
Since $\overline{U_{z_0}}^T$ is also a diagonal matrix, we have
$$\overline{U_{z_0}}^T\mathrm {diag}(\lambda_1,...,\lambda_r,0,...,0)=\mathrm {diag}(\lambda_1,...,\lambda_r,0,...,0)\overline{U_{z_0}}^T,$$
and therefore,
\begin{equation}\label{6}
\begin{split}
\overline{J_f(z_0)}^Tw_0=\mathrm {diag}(\lambda_1,...,\lambda_r,0,...,0)\overline{U_{z_0}}^T \sum_{j=1}^{r}e^n_j=\mathrm {diag}(\lambda_1,...,\lambda_r,0,...,0) z_0,
\end{split}
\end{equation}
where $\lambda_j=\frac{\partial g_1(z_0)}{\partial z_j}\geq 0$ with $\sum_{j=1}^{r}\lambda_j\geq 1$ and $1\leq j\leq r$.

{\bf Step 5.}
Let $f(z_0)=w_0$ with $z_0\in \partial D^n,~w_0\in \partial B^N$. If $f(0)=a\neq 0$, then we use the automorphism of $B^N$ to get the result. Assume $\phi_a(w)$ is an automorphism of $B^N$ such that $\phi_a(a)=0$. Then $\phi_a(w_0)\in \partial B^N$ as well. With similar analysis to Step 3, there exists a $U_{\phi_a}$ such that $U_{\phi_a}(\phi_a(w_0))=w_0$. Let $$h=U_{\phi_a}\circ \phi_a\circ f,$$ then $h(0)=0,~h(z_0)=w_0.$ As a result from Step 4, there is a sequence of real numbers $\gamma_j\geq 0$ and $\sum_{j=1}^{r}\gamma_j\geq 1$ such that
$$\overline{J_h(z_0)}^Tw_0=\mathrm {diag}(\gamma_1,...,\gamma_r,0,...,0) z_0.$$
According to the expression of $h$, it is obtained that
\begin{equation}\label{7}
\begin{split}
\overline{J_h(z_0)}^Tw_0=\overline{U_{\phi_a}J_{\phi_a}(w_0) J_f(z_0)}^Tw_0 =\overline{ J_f(z_0)}^T\overline{J_{\phi_a}(w_0) }^T\overline{U_{\phi_a}}^Tw_0.
\end{split}
\end{equation}
Since $U_{\phi_a}(\phi_a(w_0))=w_0$, then $\overline{U_{\phi_a}}^Tw_0=\phi_a(w_0)$. From the expression of the automorphism  $\phi_a$ given by \cite{rudin}, we have the following equality.
\begin{equation}\label{}
\begin{split}
\overline{J_{\phi_a}(w_0) }^T\overline{U_{\phi_a}}^Tw_0=\overline{J_{\phi_a}(w_0) }^T\phi_a(w_0)=\frac{1-||a||^2}{|1-\bar{a}^Tw_0|^2}w_0.\nonumber
\end{split}
\end{equation}
 Therefore, combining with (\ref{7}) we get
\begin{equation}\label{}
\begin{split}
\overline{ J_f(z_0)}^T\frac{1-||a||^2}{|1-\bar{a}^Tw_0|^2}w_0=\mathrm {diag}(\gamma_1,...,\gamma_r,0,...,0) z_0.\nonumber
\end{split}
\end{equation}
Consequently,
\begin{equation}\label{8}
\begin{split}
\overline{J_f(z_0)}^Tw_0=\lambda \mathrm {diag}(\gamma_1,...,\gamma_r,0,...,0) z_0,
\end{split}
\end{equation}
where $\lambda=\frac{|1-\bar{a}^Tw_0|^2}{1-||a||^2}>0$,  $\sum_{j=1}^{r}\gamma_j\geq 1$ and $ a=f(0)$.

\end{proof}

\section*{Acknowledgments}
The work was finished while the first author visited the third author at Department of Mathematical Sciences, Indiana University-Purdue University Fort Wayne. The first author would appreciate the comfortable research environment and all the support provided by the institution in academic year 2014. The work was supported by China Scholarship Council (CSC). It was also supported by the NNSF of China under Grants 11101373, and 11171255.


\end{document}